\documentclass[number,seceqn,citesort,dvips]{arxbj}
\usepackage{dcolumn}

\aid{0}
\volume{16}
\issue{4}
\pubyear{2010}
\firstpage{1164}
\lastpage{1176}
\doi{10.3150/10-BEJ257}

\makeatletter

\newtheorem{theorem}{Theorem}[section]
\newtheorem{lemma}{Lemma}[section]

\def\eqref#1{(\ref{#1})}
\newcolumntype{d}[1]{D{.}{.}{#1}}

\makeatother

\begin{document}
\begin{frontmatter}

\title{Optimal designs for discriminating between dose-response models
in toxicology studies}
\runtitle{Optimal designs for discriminating between dose-response models}

\begin{aug}
\author[1]{\fnms{Holger} \snm{Dette}\thanksref{1}\ead[label=e1]{holger.dette@ruhr-uni-bochum.de}},
\author[2]{\fnms{Andrey} \snm{Pepelyshev}\thanksref{2,e2}\ead[label=e2,mark]{andrey@ap7236.spb.edu}},
\author[2]{\fnms{Piter} \snm{Shpilev}\thanksref{2,e3}\ead[label=e3,mark]{pitsh@front.ru}}
\and\\
\author[3]{\fnms{Weng Kee} \snm{Wong}\corref{}\thanksref{3}\ead[label=e4]{wkwong@ucla.edu}}
\runauthor{Dette, Pepelyshev, Shpilev and Wong}
\address[1]{Faculty of Mathematics, Ruhr-Universit\"at Bochum,
Universitatsstrasse 150, 44780 Bochum, Germany.
\printead{e1}}
\address[2]{Department of Mathematics, St. Petersburg State University,
Universitetskij pr. 28, St. Petersburg, 198504, Russia.
E-mails: \printead*{e2,e3}}
\address[3]{Department of Biostatistics, University of California at
Los Angeles,
10833 Le Conte Avenue, Los Angeles, CA 90095, USA.
\printead{e4}}
\end{aug}

% HISTORY:
\received{\smonth{10} \syear{2008}}
\revised{\smonth{11} \syear{2009}}

% ABSTRACT
%
\begin{abstract}
We consider design issues for toxicology studies when we have a
continuous response and the true mean response is only known to be
a member of a class of nested models. This class of non-linear models was
proposed by toxicologists who were concerned only with estimation
problems. We develop robust and efficient designs for model
discrimination and for estimating parameters in the
selected model at the same time. In particular, we propose designs
that maximize the minimum of $D$- or $D_1$-efficiencies over all
models in the given class. We show that our optimal designs are
efficient for determining an appropriate model from the postulated
class, quite efficient for estimating model parameters in the
identified model and also robust with respect to model
misspecification. To facilitate the use of optimal design ideas in practice,
we have also constructed a website that freely enables practitioners to
generate a variety of optimal designs for a range of models and also
enables them to evaluate the efficiency of any design.
\end{abstract}

% KEYWORDS
%

\end{frontmatter}

%s1 ###
\section{Introduction}

This paper addresses design issues for
toxicology studies when the primary outcome is continuous and it is
not known a priori which model is an appropriate one to use. Such
design problems are common; see, for example,
\cite{1967Box,1968Hill,1975Atkinson,2007Atkinson,2009Dette}.
In
this situation, we may consider a class of plausible models within
which we believe lies an adequate model for fitting the data. The
issues of interest are how to design the study to
choose the most appropriate model from within the
postulated class of models and, at the same time, be able to estimate
the parameters of the selected model efficiently. Our design
decisions include how to select the number of dose levels to observe
the continuous outcome, where these levels are and how many repeated
observations to take at each of these levels. This work
assumes, for the sake of simplicity, that there is only one independent
variable, the dose level and only non-sequential designs are
considered.

When we have competing models, a design should be able to discriminate
among these models
and select the most appropriate ones.
Dette \cite{1990Dette,1994Dette,1995Dette} found optimal
discrimination designs for polynomial regression models, and
Dette and Roeder \cite{1997Dette} and Dette and Haller \cite{1998Dette} found
optimal discrimination designs for trigonometric and Fourier
regression models, respectively. \mbox{$T$-optimal} designs are usually used
to discriminate between homoscedastic models with normal errors
\cite{1968Hill,1967Box,1975Atkinson,2007Atkinson}. For
discriminating non-linear models,
only numerical results are possible; Lopez-Fidalgo \textit{et al.} \cite{2007Lopez}
investigated optimal designs maximizing a weighted average of two
$T$-criterion functions and Lopez-Fidalgo \textit{et al.} \cite{2008Lopez} constructed
$T$-optimal designs for
Michaelis--Menten-like models. When the design problem involves model
discrimination and another
optimality criterion, the problem is more complicated.
Hill \textit{et al.} \cite{1968Hill} was among the first to consider studies with two goals:
model discrimination and
estimation of model parameters. Dette \textit{et al.} \cite{2005Dette} gave a
concrete example
where they wanted to discriminate between the
Michaelis--Menten-model and the Emax model and estimate model
parameters in an enzyme-kinetic study. A key reason for there having
been so little research into such
design problems for non-linear models is that there are serious
technical difficulties.

The motivation for this work comes from recent proposals by
Piersma \textit{et al.} \cite{2002Piersma}, Slob \cite{2002Slob}, Slob and Pieters \cite{1998Slob}) and
Moerbeek \textit{et al.} \cite{2004Moerbeek} to use the same class of models to
study a
continuous outcome in toxicological studies. Their
interest was only in estimation problems and so they did not
consider design issues. Our purpose here is to find an optimal
design for identifying an appropriate model within the class of
models and, at the same time, provide reliable parameter estimates in the
selected model. To do this, we first find locally optimal designs
\cite{1953Chernoff}. These designs are the easiest to construct, but
they can be sensitive to nominal values and the model specification.
To overcome the risk of selecting an inappropriate model, we propose
maximin optimal designs that appear to be robust to
misspecifications in the model. These maximin optimal designs
maximize the minimum efficiency, regardless of which model in the
class of models is the appropriate one. As such, these optimal
designs provide some global protection against selecting the wrong
model from the postulated class of models. As we will show, the designs
also seem
quite robust to misspecification in the nominal values of the model
parameters.

In Section \ref{sec2}, we present some background and the proposed class of
models. We describe relationships between models in the class and
provide locally optimal designs for discriminating between plausible
models. We also show how optimal designs constructed for one set of
design parameters can be used to deduce the optimal design under
another set of design parameters. In Section~\ref{sec3}, we construct maximin
optimal designs for various subclasses of plausible models and
in Section~\ref{sec4}, we show that maximin optimal designs are robust to
misspecification of models in the postulated class. We offer a
conclusion in Section~\ref{sec5} and an \hyperref[app]{Appendix} containing technical
justifications of our results.

%s2 ###
\section{A class of dose-response models}\label{sec2}

In a general non-linear regression model, the mean response of the
outcome $Y$ is given by $ E[Y|t] = \eta(t,\theta )$, where we assume
the unknown parameter $\theta $ is $m$-dimensional. The class of
models proposed by toxicologists assumes all errors are
independent and normally distributed, and $ \eta(t,\theta )$ has one of
the following forms defined on a user-selected interval $[0,T]$:
%
%e2.5 ###
%e2.4 ###
%e2.3 ###
%e2.2 ###
%e2.1 ###
\begin{eqnarray}
\label{1.1} \eta(t,\theta) &  =  & a ;\qquad  m=1, \qquad \theta = a >0
, \\
\label{1.2}\eta(t,\theta) &  =  & a\mathrm{e}^{-bt} ;\qquad  m=2, \qquad \theta = (a,b)^T,
a>0, b>0 ,
 \\
\label{1.3} \eta(t,\theta) &  =  & a\mathrm{e}^{-bt^d} ; \qquad m=3,\qquad  \theta =
(a,b,d)^T,   a,b>0, d \geq1 ,
\\
\label{1.4} \eta(t,\theta) &  =  & a\bigl(c-(c-1)\mathrm{e}^{-bt}\bigr) ;\qquad  m=3,\qquad
\theta = (a,b,c)^T,   a,b> 0 , c \in[0,1] ,
\\
\label{1.5} \eta(t,\theta) &  =  & a\bigl(c-(c-1)\mathrm{e}^{-bt^d}\bigr) ;\qquad  m=4,\nonumber
\\
  &&\hspace*{-25pt}\theta = (a,b,c,d)^T,   a,b > 0, c \in[0,1], d\geq1   .
\end{eqnarray}

The rationale for this class of models was given for dose-response
relationships that cannot be derived from biological mechanisms. The
models are nested, in the sense that the models with a
smaller number of parameters can be obtained from another model by
setting specific values for the parameters. For instance, model
\eqref{1.5} is an extension of the models (\ref{1.4}) and
\eqref{1.3}, model~(\ref{1.3}) is an extension of model~(\ref{1.2}) and
model~\eqref{1.4} is an extension of the models~(\ref{1.2}) and
\eqref{1.1}. The hierarchy of the models is illustrated in the
following
diagram.
\begin{eqnarray*}
\begin{array}{ccccc}
(\ref{1.5})& \stackrel{d=1}{\Longrightarrow} &(\ref{1.4})
&\stackrel{c=1}{\Longrightarrow}& (\ref{1.1}) \\
\Downarrow\lefteqn{\small c=0} & & \Downarrow\lefteqn{c=0} && \\
(\ref{1.3})& \stackrel{d=1}{\Longrightarrow} &(\ref{1.2}) && \\
\end{array}
\end{eqnarray*}

We note that when $b=0$, all of the models (\ref{1.2})--(\ref{1.5}) reduce
to the constant model (\ref{1.1}), this relation not being shown in
the diagram.

Following \cite{1974Kiefer}, we consider only continuous designs.
A continuous design is
simply a probability measure $\xi$ with a finite number
of support points, say
$ t_1 , \ldots, t_n \in[0,T] $, and corresponding weights
$ \omega_1 ,\ldots, \omega_n $ with $ \omega_i > 0$ and
$\sum^n_{i=1} \omega_i = 1$.
%an experimental design.
If we fix the number of observations $N$ in advance, either by cost
or time considerations, then, roughly, $n_i=N\omega_i $ observations are
taken at point $t_i$, with $\sum^n_{i=1} n_i =N$. For many problems,
continuous designs
are easier to describe and study analytically than exact designs.

Jennrich \cite{1969Jennrich} showed that under regularity assumptions, the
asymptotic covariance matrix of the standardized least-squares
estimator $\sqrt{N/\sigma^2}$ ${\hat\theta}$ for the parameter
$\theta$ in the general non-linear model is given by the matrix
$M^{-1}(\xi, \theta )$, where
\begin{eqnarray*}%\label{2.2}
M(\xi,\theta ) = \int^T_0 f(t,\theta ) f^T (t,\theta )\,\mathrm{d}\xi(t)
\end{eqnarray*}
is the information matrix using design $\xi$ and
%
%e2.6 ###
\begin{equation}\label{f=eta/th}
f(t,\theta ) = {\partial\eta(t,\theta ) \over\partial\theta } =
( f_1(t,\theta ) ,
\ldots, f_m(t,\theta ))^T
\end{equation}
is the vector of partial derivatives of the conditional expectation
$\eta(t,\theta )$ with respect to the parameter $\theta $.
Additionally, we
consider only designs with a non-singular information matrix. A~sufficient condition for this property to hold is that the design
has $k$ support points, where $k$ is greater than or equal to the
number of
parameters in the model.

A locally optimal design maximizes a function of the information
matrix $M^{}(\xi, \theta )$ using nominal values of $\theta$ \cite{1953Chernoff}. There are several optimality criteria for estimating
purposes and for discriminating between models
\cite{1993Pukelsheim,2007Atkinson}. We are
interested in finding efficient designs for
model selection among models defined
by \eqref{1.1}--\eqref{1.5} that also provide good and robust
estimates for the parameters in the selected model. Accordingly, we construct an optimal design for pairs of
competing models that fulfills at least two of three following
requirements:
\begin{enumerate}[(3)]
\item[(1)]
The design should be able to test the hypotheses for discriminating
between two selected rival models.
For example, the hypothesis for discriminating between
the models \eqref{1.4} and \eqref{1.2}
is given by
\begin{eqnarray*}%\label{1.7}
H_0\dvtx c = 0\quad  \mbox{vs.}\quad H_1\dvtx c \in(0,1] .
\end{eqnarray*}
Van der
Vaart \cite{1998Van} described properties of test statistics
for testing such
hypotheses.
\item[(2)] The design should be able to efficiently estimate the
parameters in the corresponding pair of regression models
and for all models which are submodels of the model with the
larger number of parameters. For example, for model
(\ref{1.4}), the corresponding submodels are given by (\ref{1.2})
and (\ref{1.3}).
\item[(3)] The design should also be efficient for discriminating
between the different submodels of the model with the larger
number of parameters (which may also be nested). For example,
the optimal design for discriminating between the models
(\ref{1.2}) and (\ref{1.4}) should also be efficient for
discriminating between the models (\ref{1.1})/(\ref{1.4}) and
(\ref{1.1})/(\ref{1.2}).
\end{enumerate}

To make these ideas concrete, consider the $e_m$-optimality criterion,
where $e_m=(0,\ldots
,0,1)^T$ and $m$ is the larger of the number of parameters in the
two models under consideration. For fixed~$\theta$, a locally
$e_m$-optimal design minimizes
%e2.7 ###
\begin{equation}\label{2.5}
e_m^TM^{-1}(\xi,\theta)e_m = { \det
\widetilde{M}(\xi,\theta) \over\det M(\xi,\theta) } ,
\end{equation}
where
the matrix $M(\xi, \theta)$ is the information matrix in the model
with the larger number of parameters and the matrix $\widetilde{M}
(\xi, \theta)$ is obtained from $M (\xi, \theta)$ by deleting the
$m$th row and the $m$th column. A locally $e_m$-optimal design is just
a special case of a $c$-optimal design used for estimating
$c^T\theta$, where the vector $c$ is user-specified.

The expression \eqref{2.5}
%on the right hand side
is proportional to the asymptotic variance of the least-squares
estimate of $e_m^T \theta$, which is relevant for model
discrimination. In particular, minimizing the asymptotic variance
\eqref{2.5} provides a design with maximal power for testing a
simple hypothesis for $e_m^T \theta$. For example, if we want to
discriminate between models (\ref{1.4}) and (\ref{1.2}), we have
$m=3$ and $e_m^T \theta= (0,0,1)(a,b,c)^T=c$, and the cases $c\not=0$
and $c = 0$ give the two rival models (\ref{1.4}) and (\ref{1.2}),
respectively. Consequently, a design that minimizes the ratio in
(\ref{2.5})
is optimal for discriminating between the two models.

We next construct locally $e_m$-optimal designs for discriminating
between pairs of models~\eqref{1.3}--\eqref{1.5}.

%s2.1 ###
\subsection[Optimal discriminating designs for the models (2.2)
and (2.3)]{Optimal discriminating designs for the models \protect\eqref{1.2}
and \protect\eqref{1.3}}

For model \eqref{1.3} with $\theta= (a, b, d)^T$, the vector of
partial derivatives in \eqref{f=eta/th} is given by
%
%e2.8 ###
\begin{equation}\label{2.6}
f(t,\theta)=f(t,a,b,d)=\bigl(\mathrm{e}^{-bt^d},-at^d\mathrm{e}^{-bt^d},
-abt^d\ln(t)\mathrm{e}^{-bt^d}\bigr)^T .
\end{equation}
Our first result establishes basic
properties of locally $e_3$-optimal designs for model
\eqref{1.3}.

\begin{lemma}\label{lem2.1}
The locally $e_3$-optimal design in model \eqref{1.3}
does not
depend on the parameter $a$. Moreover, if $t_i(b, d, T)$ is a
support point of a locally $e_3$-optimal design on the interval
$[0, T]$ with corresponding weight $\omega_i(b, d, T)$, then for any
$r > 0$ and $d>0$,
\begin{eqnarray}\label{tra}
t_i(b, d, T^{1/d})&=& t_i(b, 1, T)^{1/d} ,\qquad \omega_i(b, d,
T^{1/d})= \omega_i(b, 1, T) ,\nonumber
\\[-8pt]\\[-8pt]
t_i(r b, 1, T)&=&\frac{1}{r} t_i(b, 1,r T) ,\qquad  \omega_i(r b, 1,
T)=\omega_i(b, 1,r T) .\nonumber
\end{eqnarray}
\end{lemma}

To find an efficient design for discriminating between models
\eqref{1.2} and \eqref{1.3}, we assume that the initial parameter
value of $d$ is unity. From the lemma, it is enough to calculate
locally $e_3$-optimal designs on a fixed design space for various
values of $b$ after the remaining parameters $a$ and $d$ are fixed.
Locally optimal designs on a different design space or having other
values of
the parameters can then be calculated using the relationships given in the
lemma.

To characterize locally $e_3$-optimal designs, we recall that a
set of functions
$h_1, \ldots, h_k \dvtx  I \to{\mathbb{R}}$ is
a \textit{Chebyshev system} on the interval $I$ if there exists an
$\varepsilon\in\{ -1, 1\}$ such that the inequality
\begin{eqnarray*}%\label{3.11}
\varepsilon\cdot
\left|\matrix{
h_1(t_1) &\ldots&h_1(t_k) \cr
\vdots&\ddots&\vdots\cr
h_k(t_1) &\ldots&h_k(t_k)
}\right|
> 0
\end{eqnarray*}
holds for all $t_1, \ldots, t_k \in I$ with $t_1 < t_2 < \cdots<
t_k.$ From Karlin and Studden  \cite{1966Karlin}, Theorem II 10.2, if $\{ h_1,
\ldots, h_k\}$ is a Chebyshev system, then there exists a unique
function, say $ \sum^k_{i=1} c_i^\ast h_i(t) = c^{\ast T} h(t), $
where $h=(h_1,\ldots,h_k)^T$, with the following properties:
\begin{longlist}[(ii)]
\item[(i)]
$\vert c^{\ast T}h(t)\vert\le1 \ \forall   t \in I ;$
\item[(ii)]
there exist $k$ points, $t_1^* < \cdots< t_k^*$,
such that $c^{\ast T} h(t_i^*) = (-1)^i ,  i = 1, \ldots, k.$
\end{longlist}
The function $c^{\ast T}h(t)$ alternates at the points $t^\ast_1,
\dots, t^\ast_k$ and is called the \textit{Chebyshev polynomial}.
The points
$t^*_1, \ldots, t^*_k$ are called \textit{Chebyshev points} and they are
unique when $1 \in  \operatorname{span}\{h_1, \ldots, h_k\}$, $k \ge1$ and
$I$ is a compact interval. In this case, we have $t^*_1 = \min_{t
\in I}t$, $ t^*_k = \max_{t \in I}t.$ The following result
characterizes the locally $e_3$-optimal design.

\begin{theorem}\label{thm2.1}The components of the vector defined by \eqref{2.6}
form a
Chebyshev system on the interval $[0,T]$. The locally
$e_3$-optimal design for model \eqref{1.3} is
unique and is supported at the three uniquely determined Chebyshev
points, say $t^*_1 < t^*_2 < t^*_3$.
The corresponding weights $\omega^*_1$, $\omega^*_2$, $\omega^*_3$
can be
obtained explicitly as
%
%e2.9 ###
\begin{eqnarray}
\label{weight}
\omega^* = (\omega^*_1,\omega^*_2, \omega^*_3)^T =
\frac{JF^{-1}e_3}{1_3JF^{-1}e_3} , \label{2.7}
\end{eqnarray}
where the
matrices $F$ and $J$ are defined by $ F = (f(t^*_1, \theta),
f(t^*_2, \theta), f(t^*_3, \theta))$, $J = \mathrm{diag}(1,\break -1, 1),$
respectively, and $1_3 = (1, 1, 1)^T$.
\end{theorem}

Table \ref{tab:e3-mod23} displays selected locally $e_3$-optimal
designs for model \eqref{1.3}.

%t1 ###
\begin{table}[t]
\caption{Locally $e_3$-optimal designs for models
\protect\eqref{1.3} and \protect\eqref{1.4} on
the design space $[0,1]$ for various values of the parameter $b$}\label{tab:e3-mod23}
\begin{tabular*}{\tablewidth}{@{\extracolsep{4in minus 4in}}lllllllllllll@{}}
\hline
$b $&\multicolumn{6}{l}{Model \eqref{1.3}}&\multicolumn{6}{l@{}}{Model
\eqref{1.4}}\\[-6pt]
&\multicolumn{6}{l}{\hrulefill}&\multicolumn{6}{l@{}}{\hrulefill}
\\
 & $t_1$ & $t_2$ & $t_3$ & $\omega_1$& $\omega_2$& $\omega_3$
& $t_1$ & $t_2$ & $t_3$ & $\omega_1$& $\omega_2$& $\omega_3$ \\
\hline
% $0.1$ & 0 & 0.3545 & 1 & 0.3113 & 0.4997 & 0.1890 \ \hline
% $0.2$ & 0 & 0.3415 & 1 & 0.3068 & 0.4989 & 0.1943 \ \hline
% $0.3$ & 0 & 0.3288& 1 & 0.3028& 0.4971& 0.2001 \ \hline
% $0.4$ & 0 & 0.3165& 1 & 0.2982& 0.4956& 0.2061 \ \hline
% $0.5$ & 0 & 0.3046& 1 & 0.2944& 0.4930& 0.2126 \ \hline
% $0.6$ & 0 & 0.2931& 1 & 0.2905& 0.4900& 0.2195 \ \hline
% $0.7$ & 0 & 0.2819& 1 & 0.2868& 0.4865& 0.2267 \ \hline
% $0.8$ & 0 & 0.2711& 1 & 0.2833& 0.4822& 0.2345 \ \hline
% $0.9$ & 0 & 0.2607& 1 & 0.2794& 0.4781& 0.2425 \ \hline
% $1.0$ & 0 & 0.2507& 1 & 0.2760& 0.4731& 0.2509 \ \hline
0.1&0&0.355&1 &0.311&0.500&0.189 &0&0.492&1&0.242&0.500&0.259\\
0.5&0&0.305&1 &0.294&0.493&0.213 &0&0.458&1&0.212&0.492&0.296\\
1.0&0&0.251&1 &0.276&0.473&0.251 &0&0.418&1&0.180&0.469&0.351\\
2.0&0&0.167&1 &0.241&0.403&0.356 &0&0.343&1&0.127&0.384&0.490\\
3.0&0&0.112&0.751&0.232&0.381&0.387 &0&0.281&1&0.083&0.267&0.650\\
\hline
\end{tabular*}
\end{table}

%s2.2 ###
\subsection[Optimal discriminating designs for the models (2.3) and (2.4), (2.1) and (2.4)]{Optimal
discriminating designs for the models \protect\eqref{1.3} and
\protect\eqref{1.4}, \protect\eqref{1.1} and \protect\eqref{1.4}}

For model (\ref{1.4}), we have $\theta= (a, b, c)^T$ and when $c=0$
or $c=1$, model~(\ref{1.4}) reduces to model~(\ref{1.3}) or
(\ref{1.1}), respectively. The $e_3$-optimal design is optimal for
discriminating between models (\ref{1.3}) and (\ref{1.4}) and for
discriminating between models (\ref{1.1}) and (\ref{1.4}). The
vector of partial derivatives in \eqref{f=eta/th} is
\begin{eqnarray*}%\label{2.8}
f(t,\theta)=f(t,a,b,c)=\bigl(c-(c-1)\mathrm{e}^{-bt}, a(c-1)t\mathrm{e}^{-bt},
a(1-\mathrm{e}^{-bt})\bigr)^T
\end{eqnarray*}
and its components form a Chebyshev system on $[0,T]$. The
locally $e_3$-optimal design is described in Theorem \ref{thm2.2} and we observe
that it does not depend on the parameter $a$.
For other positive values of $b$ and $T$,
the support points $t_i(b,T)$ and corresponding weights $\omega_i(b,T)$
of the optimal design
are found from $ t_i (rb, T) = \frac
{1}{r} t_i (b, rT)$ and $ \omega_i(rb, T)=\omega_i(b, rT)$.

\begin{theorem}\label{thm2.2}
Let $0\le c<1$. The locally
$e_3$-optimal design for model \eqref{1.4} on $[0,T]$ is unique and
has
three points at $t_1^*=0,$ $t_3^*=T$ and (middle point)
\begin{eqnarray*}
t_2^*=\frac{1}{b}+\frac
{t_1^*\mathrm{e}^{-bt_1^*}-t_3^*\mathrm{e}^{-bt_3^*}}{\mathrm{e}^{-bt_1^*}-\mathrm{e}^{-bt_3^*}} ,
\end{eqnarray*}
and the corresponding weights $\omega^*_1$, $\omega^*_2$ and $\omega
^*_3$ can
be obtained explicitly from formula \eqref{weight}.
% Note that we choose the points $t_1$, $t_3$ under the conditions $
%and $ \sum_{i=1}^{3}\ww_i=1$.
\end{theorem}

%the design space $[0,1]$ for various values of the parameter $b$.}
% $b $ & $t_1$ & $t_2$ & $t_3$ & $\ww_1$& $\ww_2$& $\ww_3$ \ \hline
% $0.1$ & 0 & 0.4917& 1 & 0.2418& 0.4997& 0.2585 \ \hline
% $0.2$ & 0 & 0.4834& 1 & 0.2336& 0.4992& 0.2671 \ \hline
% $0.3$ & 0 & 0.4750& 1 & 0.2263& 0.4972& 0.2765 \ \hline
% $0.4$ & 0 & 0.4668& 1 & 0.2190& 0.4950& 0.2860 \ \hline
% $0.5$ & 0 & 0.4585& 1 & 0.2119& 0.4922& 0.2959 \ \hline
% $0.6$ & 0 & 0.4503& 1 & 0.2051& 0.4888& 0.3061 \ \hline
% $0.7$ & 0 & 0.4421& 1 & 0.1985& 0.4848& 0.3167 \ \hline
% $0.8$ & 0 & 0.4340& 1 & 0.1920& 0.4802& 0.3278 \ \hline
% $0.9$ & 0 & 0.4260& 1 & 0.1858& 0.4750& 0.3392 \ \hline
% $1.0$ & 0 & 0.4180& 1 & 0.1798& 0.4692& 0.3510 \\
%0.1&0&0.492&1&0.242&0.500&0.259\\
%0.5&0&0.458&1&0.212&0.492&0.296\\
%1.0&0&0.418&1&0.180&0.469&0.351\\
%2.0&0&0.343&1&0.127&0.384&0.490\\
%3.0&0&0.281&1&0.083&0.267&0.650\\

%s2.3 ###
\subsection[Optimal discrimination designs for the models (2.4) and (2.5), (2.1) and (2.5), (2.3)
and (2.5)]{Optimal discrimination designs for the models (\protect\ref{1.4})
and (\protect\ref{1.5}), (\protect\ref{1.1}) and (\protect\ref{1.5}),
(\protect\ref{1.3}) and (\protect\ref{1.5})}

Model (\ref{1.5}) with $\theta= (a, b, c, d)^T$ reduces to model
(\ref{1.3}), (\ref{1.1}) or (\ref{1.4}) when
$c=0, c=1$ or \mbox{$ d=1$},
respectively. For testing purposes, we want an $e_3$-optimal design
for estimating the parameter $c$ and an $e_4$-optimal design for
estimating the parameter $d$. The vector of partial derivatives of
$\eta$ for
model \eqref{1.5} is
%
%e2.10 ###
\begin{equation}\label{f(tabcd)} f(t,\theta)=\bigl(c-(c-1)\mathrm{e}^{-bt^d},
a(c-1)t^d\mathrm{e}^{-bt^d},a(c-1)t^d\ln(t)b\mathrm{e}^{-bt^d},
a(1-\mathrm{e}^{-bt^d})\bigr)^T
\end{equation}
and its components form a Chebyshev
system on the interval $[0,T]$.
% In model
%(\ref{1.5}) there are two hypotheses which are of interest
%for model discrimination namely
%$H_0: c=0$ and $H_0: d=1$ under which the regression model
%reduces to the model (\ref{1.3}) and (\ref{1.4}),
%respectively. The corresponding discrimination design problems
%are given by the $e_3$- and $e_4$-optimal design problem, respectively.
Arguments similar to those given in the proof of Lemma \ref{lem2.1} show that the
support points $t_i(b, d, T)$ and weights $\omega_i(b, d, T)$ of a
locally $e_3$- or $e_4$-optimal design on the interval $[0,T]$
satisfy relations \eqref{tra}. Moreover, the optimal designs do not
depend on the parameter $a$. Table \ref{tab:e34} shows some locally
$e_3$- and $e_4$-optimal designs for model \eqref{1.5} obtained from
Theorem \ref{thm2.3} below. The proof is similar to the proof of Theorem \ref{thm2.1}
and is therefore omitted.

\begin{theorem}\label{thm2.3}The $e_3$- and $e_4$-optimal designs for model \eqref
{1.5} are
uniquely supported at the four Chebyshev points, say $t_1^* <t_2^*
< t_3^*<t_4^*$, corresponding to the Chebyshev system defined by the
components in \eqref{f(tabcd)}. The corresponding
weights, $\omega^*_1,\ldots, \omega^*_4$, are explicitly given by
\begin{eqnarray*}
\omega^*= (\omega^*_1,\ldots, \omega^*_4)^T =
\frac{JF^{-1}e_k}{1_4JF^{-1}e_k} ,\qquad   k=3,4 ,
\end{eqnarray*}
where the matrices
$F$ and $J$ are defined by $ F = (f(t^*_1, \theta), f(t^*_2,
\theta), f(t^*_3, \theta), f(t^*_4, \theta))$, $J =
\mathrm{diag}(1,-1, 1,-1),$ respectively, $1_4 = (1, 1, 1,1)^T$ and
$f(t, \theta)$ is given in \eqref{f(tabcd)}.
\end{theorem}

%t2 ###
\begin{table}[t]
\caption{Locally $e_3$- and $e_4$-optimal designs for model
\protect\eqref{1.5} on
the design space $[0,1]$ for various values of the parameter $b$}
\label{tab:e34}
\begin{tabular*}{\tablewidth}{@{\extracolsep{4in minus 4in}}lllllllllllll@{}}
\hline
$b $ & $t_1$ & $t_2$ & $t_3$ & $t_4$ &
\multicolumn{4}{l}{$e_3$-optimal}&\multicolumn{4}{l}{$e_4$-optimal}\\[-6pt]
&\multicolumn{4}{l}{}&\multicolumn{4}{l}{\hrulefill}&\multicolumn{4}{l@{}}{\hrulefill}\\
&\multicolumn{4}{l}{}& $\omega_1$& $\omega_2$&
$\omega_3$&
$\omega_4$
& $\omega_1$& $\omega_2$& $\omega_3$& $\omega_4$\\ \hline
% $0.1$ & 0 & 0.1307& 0.6475& 1 & 0.2858& 0.4156& 0.2141& 0.0844 \\
% $0.2$ & 0 & 0.1288& 0.6425& 1 & 0.2837& 0.4135& 0.2167& 0.0861 \\
% $0.3$ & 0 & 0.1268& 0.6369& 1 & 0.2811& 0.4125& 0.2189& 0.0874 \\
% $0.4$ & 0 & 0.1249& 0.6311& 1 & 0.2785& 0.4122& 0.2207& 0.0886 \\
% $0.5$ & 0 & 0.1231& 0.6257& 1 & 0.2760& 0.4112& 0.2225& 0.0903 \\
% $0.6$ & 0 & 0.1209& 0.6200& 1 & 0.2747& 0.4082& 0.2253& 0.0917 \\
% $0.7$ & 0 & 0.1189& 0.6142& 1 & 0.2729& 0.4069& 0.2271& 0.0931 \\
% $0.8$ & 0 & 0.1169& 0.6083& 1 & 0.2710& 0.4056& 0.2290& 0.0945 \\
% $0.9$ & 0 & 0.1149& 0.6021& 1 & 0.2692& 0.4042& 0.2308& 0.0958 \\
% $1.0$ & 0 & 0.1129& 0.5958& 1 & 0.2675& 0.4030& 0.2325& 0.0970 \\
0.1&0&0.131&0.648&1& 0.286&0.416&0.214&0.084& 0.174&0.328&0.326&0.172\\
0.5&0&0.123&0.626&1& 0.277&0.410&0.223&0.090& 0.156&0.302&0.342&0.200\\
1.0&0&0.113&0.596&1& 0.267&0.403&0.233&0.097& 0.137&0.272&0.352&0.239\\
2.0&0&0.094&0.530&1& 0.253&0.392&0.246&0.108& 0.106&0.215&0.341&0.338\\
3.0&0&0.079&0.463&1& 0.244&0.382&0.256&0.118& 0.080&0.163&0.289&0.468\\
\hline
\end{tabular*}
\end{table}

% $b $ & $t_1$ & $t_2$ & $t_3$ & $t_4$ & $\ww_1$& $\ww_2$& $\ww_3$& $
% $0.1$ & 0 & 0.1308& 0.6476& 1 & 0.1734& 0.3275& 0.3270& 0.1721 \\
% $0.2$ & 0 & 0.1288& 0.6423& 1 & 0.1689& 0.3215& 0.3307& 0.1789 \\
% $0.3$ & 0 & 0.1268& 0.6370& 1 & 0.1645& 0.3150& 0.3347& 0.1859 \\
% $0.4$ & 0 & 0.1247& 0.6337& 1 & 0.1605& 0.3074& 0.3382& 0.1939 \\
% $0.5$ & 0 & 0.1229& 0.6258& 1 & 0.1560& 0.3024& 0.3415& 0.2001 \\
% $0.6$ & 0 & 0.1209& 0.6200& 1 & 0.1519& 0.2962& 0.3444& 0.2075 \\
% $0.7$ & 0 & 0.1189& 0.6143& 1 & 0.1480& 0.2899& 0.3469& 0.2152 \\
% $0.8$ & 0 & 0.1169& 0.6082& 1 & 0.1442& 0.2838& 0.3489& 0.2230 \\
% $0.9$ & 0 & 0.1149& 0.6015& 1 & 0.1405& 0.2778& 0.3505& 0.2311 \\
% $1.0$ & 0 & 0.1130& 0.5957& 1 & 0.1371& 0.2717& 0.3517& 0.2395 \\
%0.1&0&0.131&0.647&1&0.174&0.328&0.326&0.172\\
%0.5&0&0.123&0.626&1&0.156&0.302&0.342&0.200\\
%1.0&0&0.113&0.596&1&0.137&0.272&0.352&0.239\\
%2.0&0&0.094&0.530&1&0.106&0.215&0.341&0.338\\
%3.0&0&0.078&0.463&1&0.080&0.163&0.289&0.468\\
%the design space $[0,1]$ for various values of the parameter $b$.}

%s3 ###
\section{Maximin optimal discriminating designs}\label{sec3}
We now wish to find an
efficient design for testing several hypotheses that
discriminate between models
(\ref{1.3}) and (\ref{1.2}), (\ref{1.4}) and (\ref{1.2}),
(\ref{1.5}) and (\ref{1.3}), and (\ref{1.5}) and (\ref{1.4}).

Let us first find an optimal design to
discriminate between two models $(2.i)$ and $(2.j)$ and let
$\operatorname{eff}^{(2.i)-(2.j)}(\xi, \theta)$ be the efficiency of the design
$\xi$ for discriminating between the two models. As an illustrative
case, consider finding the locally optimal design for
discriminating between the models (\ref{1.3}) and (\ref{1.2}). This
optimal design minimizes $e_3^TM_{(\ref{1.3})}^{-1}(\xi,\theta)e_3$
among all designs for which the matrix is regular (Theorem \ref{thm2.1}).
Here, the matrix $M_{(\ref{1.3})}(\xi, \theta)$ is the information
matrix under model (\ref{1.3}). If $\xi_3^*(\theta)$ is the locally
optimal design for discriminating between models (\ref{1.3}) and
(\ref{1.2}), then the efficiency of a design $\xi$ for discriminating
between models (\ref{1.3})--(\ref{1.2}) is defined by
\begin{eqnarray*}%\label{2.9}
\operatorname{eff}^{(\ref{1.3})-(\ref{1.2})}(\xi, \theta)=
\frac{e_3^TM_{(\ref{1.3})}^{-1}(\xi_3^*(\theta),\theta
)e_3}{e_3^TM_{(\ref{1.3})}^{-1}(\xi,
\theta)e_3} .
\end{eqnarray*}

This ratio is between $0$ and $1$; if the value is $0.5$,
this means that twice as many observations are required from the
design $\xi$ than the optimal design to discriminate between the
two models with the same level of precision. The efficiencies of
$\xi$ for discriminating between other pairs of models are
similarly defined and denoted by $
\operatorname{eff}^{(\ref{1.4})-(\ref{1.2})}(\xi, \theta)$, $
\operatorname{eff}^{(\ref{1.5})-(\ref{1.3})}(\xi, \theta)$ and $
\operatorname{eff}^{(\ref{1.5})-(\ref{1.4})}(\xi, \theta)$. Here, and
elsewhere in our work, we remind readers that we assume $\theta $ to be
fixed throughout and so all optimal designs are only
locally optimal.

Next, we use the %standardized
maximin efficient approach proposed by Dette \cite{1995Dette} and M\"uller \cite{1995Muller}
to find efficient designs for all four discrimination problems. For a
fixed $\theta $, we call a
design a \textit{maximin optimal discriminating design} for models
(\ref{1.1})--(\ref{1.5}) if it maximizes
%
%e3.1 ###
\begin{equation}\label{maximin-crit}
\min\bigl\{ \operatorname{eff}^{(\ref{1.3})-(\ref{1.2})}(\xi,
\theta), \operatorname{eff}^{(\ref{1.4})-(\ref{1.2})}(\xi, \theta),
\operatorname{eff}^{(\ref{1.5})-(\ref{1.3})}(\xi, \theta),
\operatorname{eff}^{(\ref{1.5})-(\ref{1.4})}(\xi, \theta)\bigr\}.
\end{equation}

In practice, maximin optimal discriminating designs have to be found
numerically. All computations for the optimal designs were done
sequentially using the Nelder--Mead algorithm in the {MATLAB}
package. First, maximin designs were found by maximizing the
optimality criterion within the class of all $4$-point designs. We
started with four points because that was the number of points required
for obtaining all non-zero efficiencies in \eqref{maximin-crit}.
After the
4-point optimal design was found, we searched for the optimal design
within the class of all 5-points designs and repeated the
procedure. Each time, we increased the number of points by unity,
until there was no further
improvement in the criterion value. Table \ref{tab:des-maximin} shows
maximin optimal discriminating
designs and their efficiencies when $\theta=(a,b,d,c)^T=(1,b,1,0)^T$
for different values of $b$. % for discriminating other pairs models.
We observe from the rightmost columns in the table that the maximin
optimal discriminating design has between 68--85\% efficiency
for
discriminating between %other
different
pairs of rival models from the postulated class.

%t3 ###
\begin{table}[b]
\tabcolsep=0pt
\caption{Maximin optimal discriminating designs for the
optimality criterion \protect\eqref{maximin-crit} on the design space $[0,1]$
and their efficiencies}
\label{tab:des-maximin}
\begin{tabular*}{\tablewidth}{@{\extracolsep{4in minus 4in}}lllllllllllll@{}}
\hline
$b $ & $t_1$ & $t_2$ & $t_3$ & $t_4$ & $\omega_1$& $\omega_2$&
$\omega_3$&
$\omega_4$
&    \scriptsize{\eqref{1.3}--(\ref{1.2})}    &    \scriptsize{\eqref
{1.4}--(\ref{1.2})}
&    \scriptsize{\eqref{1.5}--(\ref{1.3})}    &    \scriptsize{\eqref
{1.5}--(\ref{1.4})}
\\ \hline
% $0.1$ & 0 & 0.1815 & 0.5636 & 1 & 0.2402 & 0.2502 & 0.3201 & 0.1895 \
% \hline
% $0.2$ & 0 & 0.1878 & 0.5619 & 1 & 0.2266 & 0.2571 & 0.3220 & 0.1943 \
% \hline
% $0.3$ & 0 & 0.1885& 0.5617& 1 & 0.2292& 0.2579& 0.3114& 0.2015 \\
% $0.4$ & 0 & 0.1750& 0.5352& 1 & 0.2238& 0.2555& 0.3121& 0.2086 \\
% $0.5$ & 0 & 0.1693& 0.5326& 1 & 0.2271& 0.2545& 0.3092& 0.2093 \\
% $0.6$ & 0 & 0.1675& 0.5192& 1 & 0.2158& 0.2585& 0.3061& 0.2196 \\
% $0.7$ & 0 & 0.1683& 0.5223& 1 & 0.2075& 0.2656& 0.3000& 0.2270 \\
% $0.8$ & 0 & 0.1646& 0.5177& 1 & 0.2061& 0.2631& 0.2968& 0.2340 \\
% $0.9$ & 0 & 0.1629& 0.5173& 1 & 0.1966& 0.2613& 0.2961& 0.2461 \\
% $1.0$ & 0 & 0.1549& 0.4780& 1 & 0.2048& 0.2566& 0.2924& 0.2462 \\
0.1&0&0.175&0.552&1&0.236&0.255&0.322&0.187&0.724&0.724&0.724&0.786\\
0.5&0&0.170&0.531&1&0.220&0.260&0.308&0.212&0.719&0.719&0.719&0.787\\
1.0&0&0.160&0.507&1&0.200&0.265&0.287&0.249&0.714&0.714&0.714&0.793\\
2.0&0&0.130&0.468&1&0.161&0.250&0.249&0.340&0.705&0.702&0.702&0.848\\
3.0&0&0.105&0.440&1&0.141&0.233&0.199&0.427&0.705&0.682&0.682&0.871\\
\hline
\end{tabular*}
\end{table}

% $b $ & (\ref{1.3})-(\ref{1.2}) &
% (\ref{1.4})-(\ref{1.2}) & (\ref{1.5})-(\ref{1.3}) & (\ref{1.5})-(
% $0.1$ & 0.7230& 0.7230& 0.7230 & 0.7956 \ \hline
% $0.2$ & 0.7228& 0.7228& 0.7228 & 0.7988 \ \hline
% $0.3$ & 0.7233& 0.7179& 0.7179 & 0.7862 \ \hline
% $0.4$ & 0.7282& 0.7282& 0.7282 & 0.7750 \ \hline
% $0.5$ & 0.7197& 0.7197& 0.7197 & 0.7959 \ \hline
% $0.6$ & 0.7258& 0.7258& 0.7258 & 0.7737 \ \hline
% $0.7$ & 0.7175& 0.7175& 0.7175 &0.7871 \ \hline
% $0.8$ & 0.7163& 0.7163& 0.7163 &0.7917 \ \hline
% $0.9$ & 0.7409& 0.7159& 0.7159 & 0.8026 \ \hline
% $1.0$ & 0.7354& 0.7354& 0.7354 &0.7450 \ \hline
%designs in Table 5 for discriminating %other
%different
%pairs of models.}

%s4 ###
\section{Efficiencies of maximin optimal designs for estimating model
parameters under model uncertainty}\label{sec4}

We now investigate the performance of maximin discrimination designs
for estimating parameters in the different models. We first present
results for estimating each parameter in the model and
$D$-efficiencies of the maximin discrimination design for estimating
all parameters in the model. We recall that $D$-efficiencies are
computed relative to the $D$-optimal design for the specific model
and $D$-optimal designs are found by maximizing the determinant of
the expected information matrix over all designs on the design
space. $D$-optimal designs are appealing because they minimize the
generalized variance and thereby provide the smallest volume of the confidence
ellipsoid for all parameters in the mean function.

Table \ref{tabind-eff} displays efficiencies of selected maximin
optimal discriminating designs for estimating the individual
parameters in the four models. The efficiencies for
estimating the parameter $a$ are consistently the lowest and
efficiencies for estimating the parameters $b$, $c$
and $d$ tend to be sequentially higher for each model. It is
not surprising to observe that the efficiencies are highest
for estimating the particular parameter that sets the two models
apart.

%t4 ###
\begin{table}[b]
 \caption{Efficiencies of the maximin optimal
discriminating designs in Table \protect\ref{tab:des-maximin} for estimating
individual coefficients in models \protect\eqref{1.2}--\protect\eqref{1.5}; the
first two columns are efficiencies
for estimating $a$ and $b$ in model~\protect\eqref{1.2}, the next three
columns are for estimating $a$, $b$ and $d$ in model~\protect\eqref{1.3}, the
next three columns are efficiencies for estimating $a$, $b$ and $c$
in model \protect\eqref{1.4} and the last four columns are for estimating
$a$, $b$, $c$ and $d$ in model~\protect\eqref{1.5}}\label{tabind-eff}
\begin{tabular*}{\tablewidth}{@{\extracolsep{4in minus 4in}}lllllllllllll@{}}
\hline
$b $ &\multicolumn{2}{l}{Model \eqref{1.2}}&\multicolumn{3}{l}{Model \eqref{1.3}}&
\multicolumn{3}{l}{Model \eqref{1.4}}&\multicolumn{4}{l}{Model \eqref{1.5}}
\\[-6pt]
&\multicolumn{2}{l}{\hrulefill}&\multicolumn{3}{l}{\hrulefill}&
\multicolumn{3}{l}{\hrulefill}&\multicolumn{4}{l@{}}{\hrulefill}
\\
& $\operatorname{eff}_1$ & $\operatorname{eff}_2$ & $\operatorname{eff}_1$ & $\operatorname{eff}_2$ & $\operatorname{eff}_3$
& $\operatorname{eff}_1$ & $\operatorname{eff}_2$ & $\operatorname{eff}_3$ & $\operatorname{eff}_1$ &
 $\operatorname{eff}_2$ & $\operatorname{eff}_3$& $\operatorname{eff}_4$ \\
\hline
% $0.1 $ & 0.4188 & 0.4982 & 0.2642 & 0.4980 & 0.7230 \ \hline
% $0.2 $ & 0.4002 & 0.4918 & 0.2486 & 0.4913 & 0.7228 \ \hline
% $0.3 $ & 0.3947 & 0.5078 & 0.2551 & 0.5043 & 0.7233 \ \hline
% $0.4 $ & 0.3759 & 0.4986 & 0.2406 & 0.4929 & 0.7282 \ \hline
% $0.5 $ & 0.3725 & 0.5057 & 0.2431 & 0.4932 & 0.7197 \ \hline
% $0.6 $ & 0.3576 & 0.5039 & 0.2315 & 0.4878 & 0.7258 \ \hline
% $0.7 $ & 0.3439 & 0.5115 & 0.2231 & 0.4892 & 0.7175 \ \hline
% $0.8 $ & 0.3359 & 0.5222 & 0.2198 & 0.4919 & 0.7163 \ \hline
% $0.9 $ & 0.3113 & 0.5387 & 0.2006 & 0.5040 & 0.7409 \ \hline
% $1.0 $ & 0.3275 & 0.5346 & 0.2195 & 0.4828 & 0.7354 \ \hline
0.1&0.42&0.495&0.26&0.495&0.724&0.30&0.726&0.724&0.24&0.784&0.724&0.786\\
0.5&0.37&0.501&0.25&0.490&0.719&0.27&0.725&0.719&0.22&0.773&0.719&0.787\\
1.0&0.32&0.545&0.22&0.495&0.714&0.25&0.716&0.714&0.20&0.760&0.714&0.793\\
2.0&0.24&0.609&0.17&0.325&0.705&0.21&0.661&0.702&0.16&0.759&0.702&0.849\\
3.0&0.20&0.429&0.15&0.514&0.705&0.18&0.560&0.682&0.14&0.708&0.682&0.871\\
\hline
\end{tabular*}
\end{table}

% $b $ & $\operatorname{eff}_1(2.5)$ & $\operatorname{eff}_2(2.5)$ & $\operatorname{eff}_3(2.5)$& $\operatorname{eff}_1(2.6)$ &
%$\operatorname{eff}_2(2.6)$ & $\operatorname{eff}_3(2.6)$ & $\operatorname{eff}_4(2.6)$ \ %\hline
% $0.1 $ & 0.3024 & 0.7258 & 0.7230 & 0.2432 & 0.7838 & 0.7230 & 0.7956
% $0.2 $ & 0.2866 & 0.7276 & 0.7228 & 0.2288 & 0.7779 & 0.7228 & 0.7988
% $0.3 $ & 0.2919 & 0.7257 & 0.7179 & 0.2371 & 0.7782 & 0.7179 & 0.7862
% $0.4 $ & 0.2761 & 0.7330 & 0.7282 & 0.2238 & 0.7598 & 0.7282 & 0.7750
% $0.5 $ & 0.2796 & 0.7267 & 0.7197 & 0.2267 & 0.7711 & 0.7197 & 0.7959
% $0.6 $ & 0.2675 & 0.7314 & 0.7258 & 0.2158 & 0.7563 & 0.7258 & 0.7737
% $0.7 $ & 0.2588 & 0.7226 & 0.7175 & 0.2089 & 0.7680 & 0.7175 & 0.7871
% $0.8 $ & 0.2557 & 0.7210 & 0.7163 & 0.2060 & 0.7685 & 0.7163 & 0.7917
% $0.9 $ & 0.2373 & 0.7172 & 0.7159 & 0.1869 & 0.7967 & 0.7159 & 0.8026
% $1.0$ & 0.2541 & 0.7362 & 0.7354 & 0.2059 & 0.7104 & 0.7354 & 0.7450 \
% \hline
%discriminating designs in Table 5 for estimating the individual
%coefficients in model (\ref{1.4}) and (\ref{1.5}). The first three
%columns are efficiencies for estimating $a$, $b$ and $c$ in model (2.5)
%and the last four columns are for estimating $a$, $b$, $c$ and $d$ in
%model (2.6).}

Table \ref{tab:D-eff} shows $D$-efficiencies of the maximin optimal
discriminating designs in Table \ref{tab:des-maximin}. These are
efficiencies relative to each of the locally $D$-optimal designs
found for each model in the class. For the values of $b$ in Table
\ref{tab:D-eff}, all efficiencies are high. Recall that the optimal
discriminating designs were constructed for discriminating between
models \eqref{1.2} and \eqref{1.3}. We observe that these
efficiencies are highest for the most complicated model, \eqref{1.5},
averaging 96\%, while the efficiencies are about 67\% for the least
complicated model, \eqref{1.2}. This implies that the maximin
optimal designs are quite robust to misspecification of models
within the class of models and also quite insensitive to small
changes to the nominal values of the parameter $b$ common to all of the
models. In models \eqref{1.2} and \eqref{1.3}, the $D$-efficiencies
drop by roughly 15\% when the nominal value of $b$ is increased
from $2$ to $3$.

%t5 ###
\begin{table}[t]
\caption{$D$-efficiencies of maximin designs in Table
\protect\ref{tab:des-maximin}
under %other
various model assumptions} \label{tab:D-eff}
\begin{tabular*}{200pt}{@{\extracolsep{4in minus 4in}}lllll@{}}
\hline
$b $ & $\operatorname{eff}_D^{\eqref{1.2}}$ & $\operatorname{eff}_D^{\eqref{1.3}}$ &
$\operatorname{eff}_D^{\eqref{1.4}}$ & $\operatorname{eff}_D^{\eqref{1.5}}$
\\ \hline
% 0.1 & 0.7123 & 0.8521 & 0.8521 & 0.9638 \ \hline
%0.2 & 0.7136 & 0.8525 & 0.8524 & 0.9623 \ \hline
%0.3 & 0.7307 & 0.8630 & 0.8627 & 0.9711 \ \hline
%0.4 & 0.7289 & 0.8622 & 0.8617 & 0.9649 \ \hline
%0.5 & 0.7404 & 0.8650 & 0.8641 & 0.9706 \ \hline
%0.6 & 0.7433 & 0.8657 & 0.8648 & 0.9664 \ \hline
%0.7 & 0.7514 & 0.8649 & 0.8630 & 0.9686 \ \hline
%0.8 & 0.7629 & 0.8676 & 0.8657 & 0.9706 \ \hline
%0.9 & 0.7685 & 0.8570 & 0.8637 & 0.9686 \ \hline
%1.0 & 0.7846 & 0.8804 & 0.8762 & 0.9627 \ \hline
0.1&0.710&0.851&0.851&0.963\\
0.5&0.737&0.862&0.861&0.968\\
1.0&0.786&0.873&0.869&0.972\\
2.0&0.703&0.864&0.860&0.959\\
3.0&0.525&0.716&0.820&0.917\\
\hline
\end{tabular*}
\end{table}

%s5 ###
\section{Conclusions}\label{sec5}

Our work is motivated by toxicologists' recent interest in a
class of non-linear nested models for studying a continuous outcome.
The toxicologists were primarily interested in estimating parameters
or a function of model parameters. The designs employed in their
studies lacked justification. Our work addresses design issues for
such a problem, where there is model uncertainty and all candidate
models are non-linear models nested within one another. The proposed
optimal designs are efficient for model discrimination and
parameter estimation. Previous design work for discriminating
between non-linear models usually focused on two rival models; our
work finds efficient and analytic locally optimal discriminating
designs for discriminating between
pairs of models within the predetermined class.

Our proposed optimal designs were constructed using large-sample
theory. The variances of the estimated parameters were obtained via the
asymptotic covariance matrix that our optimal designs used to minimize the
asymptotic variances. It is reasonable to ask whether the asymptotic
variance is a good approximation to the actual variance of the
estimated parameters encountered in practice with realistic sample
size. We performed a small simulation study using the setup in~\cite{2002Piersma}, where rats were prenatally exposed to
diethylstilbestrol and the design $\xi_u$ had 6 animals in each of
the 10 dose groups at 0, 1.0, 1.7, 2.8, 4.7, 7.8, 13, 22, 36 and 60
mg$/$kg body weight per day. In total, there were 60 observations
from the dose interval $[0,60]$. The maximin optimal design
$\xi_{mm}$ for $b=0.1$, $d=1$, $c=0$ requires 7 rats at the $0$ dose,
$12$ rats at the $3.6$ dose, $13$ rats at the $24$ dose and $28$ at
the $60$ dose.

We simulated data with $a=1$, $\sigma=0.05$ and several values of
the parameters $b,$ $d$ and $c$. A~total of 1000 repetitions were used
in each simulation. In Table \ref{tab:simvar}, we report simulated
normalized variances of least-squares estimated parameters that are
most important for discrimination. We see that in all of the cases we
investigated, the variances using the maximin optimal design
$\xi_{mm}$ are smaller than the variances obtained from the design
$\xi_u$ of Piersma \textit{et
al.}, in many cases by a huge margin. This shows
the benefits of incorporating optimal design ideas into the design of
a toxicology study. The design of Piersma \textit{et al.} was not theory-based
and required more dose levels, which usually translated to higher
labor, material and time costs without gain in precision for the
estimates relative to the optimal design. Additional simulation
results not shown here confirm that the asymptotic variances
are close to the simulated variances.

Finally, we mention that, in principle, the approach presented here can
be applied to discriminate between models when the difference between the
dimensions of the null hypothesis and the alternative is greater than
$1$. For example, suppose that we wish to discriminate between model
(\ref{1.2}) and
the model with mean response given by $a + b_1\mathrm{e}^{c_1t} + b_2t$. In
this case,
the design maximizing the non-centrality parameter of the likelihood
ratio test
depends on the values of the parameters of the larger model. The extension
of our procedure to $D_s$-optimal designs for minimizing the volume
of the confidence ellipsoid for the parameters $(b_1,b_2)$ would
still work,
but some efficiency would be lost.

%t6 ###
\begin{table}[t]
\caption{Simulated normalized variances of some parameters in
models (\protect\ref{1.4})--(\protect\ref{f=eta/th}) for several true values of parameters (left three
columns)} \label{tab:simvar}
\begin{tabular*}{\tablewidth}{@{\extracolsep{4in minus 4in}}lllllllllll@{}}
\hline
$b $ & $d$ & $c$& \multicolumn{4}{l}{Maximin design $\xi_{mm}$}&\multicolumn{4}{l@{}}{Design $\xi_u$}\\[-6pt]
\multicolumn{3}{l}{}& \multicolumn{4}{l}{\hrulefill}&\multicolumn{4}{l@{}}{\hrulefill}\\
\multicolumn{3}{l}{}&{(\ref{1.4})}&{(\ref{1.5})}&{(\ref{f=eta/th})}&{(\ref{f=eta/th})}
&{(\ref{1.4})}&{(\ref{1.5})}&{(\ref{f=eta/th})}&{(\ref{f=eta/th})}\\
\multicolumn{3}{l}{}
& \multicolumn{1}{l}{$\mathrm{var} (\hat d)$}& \multicolumn{1}{l}{$\mathrm{var} (\hat c)$}
& \multicolumn{1}{l}{$\mathrm{var} (\hat d)$}& \multicolumn{1}{l}{$\mathrm{var} (\hat c)$}
& \multicolumn{1}{l}{$\mathrm{var} (\hat d)$}& \multicolumn{1}{l}{$\mathrm{var} (\hat c)$}
& \multicolumn{1}{l}{$\mathrm{var} (\hat d)$}& \multicolumn{1}{l@{}}{$\mathrm{var} (\hat c)$}
 \\ \hline
0.10&1.0&0.0&58.85& 2.02& \phantom{0}71.07& \phantom{0}2.48&62.73& \phantom{0}5.53& \phantom{0}88.42& \phantom{0}7.81\\
0.10&0.8&0.0&30.38& 5.39& \phantom{0}83.93&28.91&34.93&26.72& \phantom{0}86.34&58.80\\
0.10&1.0&0.2&11.86& 1.96&103.40& \phantom{0}2.79&18.43& \phantom{0}4.83&138.97& \phantom{0}7.77\\
0.10&0.8&0.2&19.39& 4.21&135.00&35.91&23.57&10.88&148.29&81.56\\
0.06&1.0&0.0&61.67& 3.91&103.78& \phantom{0}7.17&61.62&11.35&115.70&23.33\\
0.08&1.0&0.1&22.58& 2.47& \phantom{0}92.44& \phantom{0}3.68&36.83& \phantom{0}6.48&113.39&10.80\\
\hline
\end{tabular*}
\end{table}

To facilitate the use of optimal designs for practitioners, we have
created a website that freely generates different types of
tailor-made optimal designs for various popular models.

We are currently
refining the computer algorithms for generating optimal designs
discussed here and plan to upload them to the site at
\href{http://optimal-design.biostat.ucla.edu/optimal}{http://optimal-design.biostat.ucla.edu/}
\href{http://optimal-design.biostat.ucla.edu/optimal}{optimal}. We hope that the
site will stimulate interest in design issues, inform practitioners
and enable them to incorporate optimal design ideas into their work.

\begin{appendix}
\section*{\texorpdfstring{Appendix: Proofs of Lemma \protect\ref{lem2.1} and
Theorem \protect\ref{thm2.1}}{Appendix: Proofs of Lemma 2.1 and Theorem 2.1}}\label{app}

%s5.1 ###
\subsection{\texorpdfstring{Proof of Lemma \protect\ref{lem2.1}}{Proof of Lemma 2.1}}
Let $I(t, a, b, d) = f(t, a, b, d)f^T(t, a, b, d)$, where $f(t,a,b,d)$
is given in \eqref{2.6}. Lemma \ref{lem2.1} follows
from the identities
\begin{eqnarray*}
\det\int_{0}^{T}I(t, a,
b, d)\,\mathrm{d}\xi(t)=\gamma\det\int_{0}^{T^d}I(t^d, 1, b,
1)\,\mathrm{d}\xi(t)=\gamma\det\int_{0}^{T}I(t, 1, b, 1)\,\mathrm{d}\xi(t^{1/d})
\end{eqnarray*}
and
\begin{eqnarray*}\det\int_{0}^{T}I(t, a, rb,
1)\,\mathrm{d}\xi(t)=\gamma' \det\int_{0}^{T}I(rt, 1, b, 1)\,\mathrm{d}\xi(t)=\gamma'
\det
\int_{0}^{T}I(t, 1, b, 1)\,\mathrm{d}\xi(t/r) ,
\end{eqnarray*}
where $\gamma$ and
$\gamma'$ denote
appropriate constants.

%s5.2 ###
\subsection{\texorpdfstring{Proof of Theorem \protect\ref{thm2.1}}{Proof of Theorem 2.1}}
Let $g(t) = p^T f(t)$ be an arbitrary linear combination of the
functions $\mathrm{e}^{-bt},-t\mathrm{e}^{-bt}$ and $-t\ln(t)\mathrm{e}^{-bt}$. One can show that
$(g(t)\mathrm{e}^{bt})'' = c/t$ does not have any roots in the interval $[0,
T]$ and so the function $g(t)$ has at most two roots. This proves that
the system of functions $\mathrm{e}^{-bt},-t\mathrm{e}^{-bt},-t\ln(t)\mathrm{e}^{-bt}$ has the
Chebyshev property and this argument also shows that there exist
precisely three Chebyshev
points.

The proof of the remaining part now follows by a standard argument
in classical optimal design theory. After showing that the functions
$\mathrm{e}^{-bt},-t\mathrm{e}^{-bt}$ form a Chebyshev system on the interval $[0,
T]$, we have
\[
\left|\matrix{
\mathrm{e}^{-bt_1} & \mathrm{e}^{-bt_2} & 0 \cr
-t_1\mathrm{e}^{-b t_1} & -t_2\mathrm{e}^{-b t_2} & 0 \cr
-t_1\ln(t_1)\mathrm{e}^{-b t_1} & -t_2\ln(t_2)\mathrm{e}^{-bt_2} & 1\cr
}\right|
=
\left|\matrix{
\mathrm{e}^{-bt_1} & \mathrm{e}^{-bt_2} \cr
-t_1\mathrm{e}^{-bt_1} & -t_2\mathrm{e}^{-bt_2} \cr
}\right|\neq 0
\]
for all $0\leq t_1 < t_2 \leq T$ and, consequently, from
\cite{1966Karlin}, Theorem 7.7, the locally $e_3$-optimal
design is supported at the Chebyshev points. The assertion on the
weights of the locally $e_3$-optimal design follows from \cite{1991Pukelsheim}.

\end{appendix}

\section*{Acknowledgements} The authors are grateful to Martina Stein,
who typed parts of this paper with considerable technical expertise.
The work of H. Dette was supported by the Collaborative
Research Center ``Statistical modeling of non-linear dynamic processes''
(SFB 823) of the German Research Foundation (DFG) and in part by a NIH
Grant award IR01GM072876 and the BMBF-Grant SKAVOE. The work of
W.K. Wong was partially supported by NIH Grant awards R01GM072876, P01
CA109091 and P30 CA16042-33. The work of A. Pepelyshev and P. Shpilev
was partially
supported by an RFBR Grant, project number 09-01-00508.

\printhistory

\end{document}